\documentclass[12pt]{amsproc}

\usepackage[top=1in,bottom=1in,left=1.5in,right=1.5in]{geometry}
\usepackage{amsfonts} 
\usepackage{amsmath}
\usepackage{amssymb}
\usepackage{amsthm}
\usepackage{mathrsfs}
\usepackage{setspace}
\usepackage{subfigure}
\usepackage{url}
\usepackage{booktabs}
\usepackage{ifthen}
\usepackage{tikz}
\usetikzlibrary{calc}
\usepackage{enumerate}
\usepackage{enumitem}

\newcommand{\C}{\mathscr{C}}

\newcommand{\cart}{\,\square \,}

\newtheorem{thm}{Theorem}[section]

\newtheorem{prop}[thm]{Proposition}
\newtheorem{cor}[thm]{Corollary}

\newtheorem{claim}[thm]{Claim}

\newtheorem{defn}[thm]{Definition}

\numberwithin{equation}{section}

\begin{document}

\title[Vizing's conjecture]{A new bound for Vizing's conjecture}

\author{Elliot Krop \and Kimber Wolff}
\address[Elliot Krop]{Department of Mathematics, Clayton State University, Morrow, Georgia, USA 302620}
\address[Kimber Wolff]{Department of Mathematics, Vanderbilt University, Nashville, Tennessee, USA 37240}
\email[Elliot Krop]{elliotkrop@clayton.edu}
\email[Kimber Wolff]{kimber.wolff@vanderbilt.edu}
\date{\today}

\begin {abstract}
For any graph $G$, we define the power $\pi(G)$ as the minimum of the largest number of neighbors in a $\gamma$-set of $G$, of any vertex, taken over all $\gamma$-sets of $G$. We show that $\gamma(G\square H)\geq \frac{\pi(G)}{2\pi(G) -1}\gamma(G)\gamma(H)$.\\ 
Our methods allow us to prove the following statements for any graphs $G$ and $H$, 
\begin{enumerate} 
\item $\gamma(G\cart H)\geq \frac{\lceil \frac{\gamma (G)}{2}\rceil}{2\lceil \frac{\gamma (G)}{2}\rceil-1}\gamma(G)\gamma(H)$ for odd $\gamma(G)$,
\item $\gamma(G\cart H)\geq \frac{\gamma (G)}{2\gamma (G)-2}\gamma(G)\gamma(H)$, for even $\gamma(G)$, and
\item a short proof of Vizing's conjecture where $\gamma(G)=3$.
\end{enumerate}
\smallskip
Our argument relies on establishing efficient correspondences between dominating vertices and subsets of their neighborhoods and then showing a sufficient number of dominating vertices that horizontally dominate vertically undominated cells. 
\\[\baselineskip] 
\\[\baselineskip]
\end {abstract}

\maketitle

 \section{Introduction}
The famous conjecture of Vadim G. Vizing (1963) \cite{Vizing, V1} is the natural statement that for any two graphs $G$ and $H$,
\begin{align}
\gamma(G \cart H) \geq \gamma(G)\gamma(H).\label{V}
\end{align}

The survey \cite{BDGHHKR} discusses many past results and contemporary approaches to the problem. For more recent partial results see \cite{Wolff}, \cite{Zerbib}, \cite{B}, \cite{B2}, \cite{CK}, \cite{PPS}, \cite{BH}, and \cite{K}.

We say a bound is of \emph{Vizing-type} with constant $c$ if for all graphs $G$ and $H$, and a fixed constant $c$, it is shown that $\gamma(G\cart H)\ge c \gamma(G)\gamma(H)$. The first Vizing-type bound was shown in 2000 by Clark and Suen \cite{CS} with a constant $\frac{1}{2}$. This bound was improved by Zerbib in 2019 \cite{Zerbib} to 
\begin{align}
\gamma(G\cart H)\ge \frac{1}{2}\gamma(G)\gamma(H) +\frac{1}{2}\max\{\gamma(G), \gamma(H)\}. \label{Zerbibbound}
\end{align}

The best current general bound for the conjectured inequality was shown in 2020 by Wolff \cite{Wolff}. In most cases, it combined the arguments of Bre\v{s}ar \cite{B2} and Zerbib \cite{Zerbib} and improved on their previous result. Wolff showed that for any graphs $G$ and $H$, if $\rho(H)$ is the two-packing number of $H$,
\begin{align}
\gamma(G\cart H)\geq \gamma(G)\frac{2\gamma(H)-\rho(H)+1}{3}. \label{Wolff}
\end{align}

\medskip

For specific families of graphs, some better bounds exist. Krop \cite{K} showed that if $G$ is a claw-free graph and $H$ is any graph, then

\[\gamma(G\cart H)\ge \frac{2}{3}\gamma(G)\gamma(H).\]

This was improved by Bre\v{s}ar and Henning \cite{BH} to

\[\gamma(G\cart H)\ge \frac{3}{4}\gamma(G)\gamma(H).\]

\medskip

Moreover, Bre\v{s}ar \cite{B} has shown that Vizing's conjecture holds for any graphs $G$ and $H$ where $\gamma(G)=3$.

\medskip

In this paper, we improve and generalize the methods of \cite{K}. We define the \emph{power} of a graph $\pi(G)$ which is bounded above by $\gamma(G)$, and generalize the techniques of \cite{K} to show that for any graphs $G$ and $H$, 
\begin{align}
\gamma(G\square H)\geq \frac{\pi(G)}{2\pi(G)-1}\gamma(G)\gamma(H).
\end{align}

Our argument allows us to show
\begin{enumerate} 
\item $\gamma(G\cart H)\geq \frac{\lceil \frac{\gamma (G)}{2}\rceil}{2\lceil \frac{\gamma (G)}{2}\rceil-1}\gamma(G)\gamma(H)$ for odd $\gamma(G)$,
\item $\gamma(G\cart H)\geq \frac{\gamma (G)}{2\gamma (G)-2}\gamma(G)\gamma(H)$, for even $\gamma(G)$, and
\item a short proof of Vizing's conjecture where $\gamma(G)=3$.
\end{enumerate}
\medskip
Our argument relies on establishing efficient correspondences between dominating vertices and subsets of their neighborhoods and then showing a sufficient number of dominating vertices that horizontally dominate vertically undominated cells. 

\subsection{Notation and past results}

\subsubsection{Basic Definitions}

\bigskip

The graphs we discuss are finite, simple, connected, undirected graphs with vertex set $V$ and edge set $E$. We may refer to the vertex set and edge set of $G$ as $V(G)$ and $E(G)$, respectively.  For more on basic graph theoretic notation and definitions, we refer to Diestel~\cite{Diest}. 
 
For any graph $G(V,E)$, a subset $S\subseteq V$ \emph{dominates} $G$ if $N[S]=V$. The minimum cardinality of $S \subseteq V$, so that $S$ dominates $G$ is called the \emph{domination number} of $G$ and is denoted $\gamma(G)$. We call a dominating set that realizes the domination number a $\gamma$-set.

 The \emph{Cartesian product} of two graphs $G_1(V_1,E_1)$ and $G_2(V_2,E_2)$, denoted by $G_1 \cart G_2$, is a graph with vertex set $V_1 \times V_2$ and edge set $E(G_1 \cart G_2) = \{((u_1,v_1),(u_2,v_2)) : v_1=v_2 \mbox{ and } (u_1,u_2) \in E_1, \mbox{ or } u_1 = u_2 \mbox{ and } (v_1,v_2) \in E_2\}$.


Let $\Gamma=\{v_1,\dots, v_k\}$ be a minimum dominating set of $G$. For any $i\in [k]$, define the set of \emph{private neighbors} for $v_i$, $P_i=\big\{v\in V(G)-\Gamma: N(v)\cap \Gamma = \{v_i\}\big\}$. For $S\subseteq [k]$, $|S|\geq 2$, we define the \emph{shared neighbors} of $\{v_i:i\in S\}$, $P_S=\big\{v\in V(G)-\Gamma: N(v)\cap \Gamma=\{v_i: i\in S\}\big\}$.

For any $S\subseteq [k]$, say $S=\{i_1,\dots, i_s\}$ where $s\geq 2$, we may write $P_S$ as $P_{\{i_1,\dots, i_s\}}$ or $P_{i_1,\dots, i_s}$ interchangeably.

\bigskip

\subsubsection{Cells, chambers, and vertical domination}

For $i\in [k]$, let $Q_i=\{v_i\} \cup P_i$. We call $\mathcal{Q}=\{Q_1,\dots, Q_k\}$ the \emph{cells} of $G$. For any $I\subseteq [k]$, we write $Q_I=\bigcup_{i\in I}Q_i$ and call $\C(\cup_{i\in I}Q_i)=\bigcup_{i\in I}Q_i\cup\bigcup_{S\subseteq I}P_{S}$ the \emph{chamber} of $Q_I$. We may write this as $\C_{I}$.

In Figure \ref{chamber} below, the black vertices are in the minimum dominating set. The chamber of $Q_{1,2,3}$ is composed of the black and pink vertices.

\begin{figure}[ht]
\begin{center}
\begin{tikzpicture}[scale=1.5]
\tikzstyle{vert}=[circle,fill=black,inner sep=3pt]
\tikzstyle{overt}=[circle, draw, inner sep=3pt, minimum size=6pt]
\tikzstyle{rvert}=[circle,draw,fill=red!30,inner sep=3pt]

\node[vert, label=below:\tiny{$v_1$}] (v1) at (1,2) {};
\node[vert, label=below:\tiny{$v_2$}] (v2) at (2,2) {}; 
\node[vert, label=below:\tiny{$v_3$}] (v3) at (3,2) {};
\node[overt, label=below:\tiny{$v_4$}] (v4) at (4,2) {};

\node[rvert] (u1) at (.5,1.5) {};
\node[rvert] (u2) at (1.5,1.5) {};
\node[rvert] (u3) at (1.5,1) {};
\node[rvert] (u4) at (2.5,1.5) {};
\node[rvert] (u5) at (3.5,1) {};
\node[overt] (u6) at (3.5,1.5) {};
\node[overt] (u7) at (4.5,1.5) {};
\node[overt] (u8) at (3,1) {};

\draw[color=black] 
   (u1)--(v1)--(u2)--(v2) (v1)--(u3)--(v2) (v2)--(u4)--(v3) (v3)--(u5) (v3)--(u6)--(v4) (v4)--(u7) (v2).. controls (2.5,1.75) .. (u8)..controls (3.5,1.75)..(v4);

\end{tikzpicture}
\caption{The Chamber of $Q_{1,2,3}$}
\label{chamber}
\end{center}
\end {figure}

For a vertex $h\in V(H)$, the \emph{$G$-fiber}, $G^h$, is the subgraph of $G\cart H$ induced by $\{(g,h):g\in V(G)\}$. 

Since for any $h\in V(H)$, $G^h$ is isomorphic to $G$, we may refer to the cells of $G^h$, meaning those sets of vertices in $G^h$, which when projected onto $G$ correspond to cells of $G$.

For a minimum dominating set $D$ of $G\cart H$, we define $D^h=D\cap V(G^h)$. Likewise, for any set $S\subseteq [k]$, $P_S^h=P_S \times \{h\}$, and for $i\in [k]$, $Q_i^h=Q_i\times \{h\}$. For any cell $Q_i^h$ the \emph{index} is defined to be $i$. By $v_i^h$ we mean the vertex $(v_i,h)$. For any $I^h\subseteq [k]$, where $I^h$ represents the indices of some cells in $G$-fiber $G^h$, we write $\C_{I^h}$ to mean the chamber of $Q^h_{I^h}$, that is, the set $\bigcup_{i\in I^h}Q^h_i\cup\bigcup_{S\subseteq I^h}P^h_{S}$.

For clarity, assume that our representation of $G\cart H$ is with $G$ on the $x$-axis and $H$ on the $y$-axis.

Any vertex $v=(g,h)\in V(G)\times V(H)$ is \emph{vertically dominated} if $(\{g\}\times N_H[h])\cap D \neq \emptyset$ and \emph{vertically undominated}, otherwise. For $i\in [k]$ and $h\in V(H)$, we say that the cell $Q_i^h$ is \emph{vertically dominated} if $(Q_i\times N_H[h])\cap D\neq\emptyset$. A cell which is not vertically dominated is \emph{vertically undominated}. 

\bigskip





\subsubsection{The power of a graph}

\begin{defn}
For a fixed $\gamma$-set $D$ of $G$, we call $a_G(D)=\max_{v\in V(G)}\{|D\cap N[v]|\}$ the \emph{allegiance} of $D$ with respect to $G$.
\end{defn}

\begin{defn}
We call the $\pi(G)=\min_{D}\{a_G(D)\}$ taken over all $\gamma$-sets $D$ of $G$, the \emph{power} of a graph $G$.
\end{defn}

This concept is related to $[1,k]$-sets defined by Chellali et al. \cite{CHHM}.

\medskip

As an example of the power of a graph, notice that for any claw-free graph $G$, a theorem of Alan and Lasker \cite{AL} states that there exists a minimum dominating set of $G$, $D$, whose vertices are independent. Since $G$ is claw-free, $a_G(D)\le 2$ which means $\pi(G)\le 2$.

\medskip

%
%

\subsection{A Useful Inequality}
Although the following inequality is elementary, we provide the proof for completeness.

\begin{prop}\label{ineq}
If 
\begin{align*}
f(t_1,\dots, t_n)=\sum_{i=1}^n{i\times t_i}
\end{align*}
subject to 
\begin{align}
&\sum_{i=1}^n{t_i}=1 \label{i}\\
&t_1\geq \sum_{i=2}^n{(i-1)t_i}, \label{ii}
\end{align} 
for real-valued $t_i, 1\leq i \leq n$, then 
\begin{align*}f(t_1,\dots, t_n)\leq \frac{2n-1}{n}
\end{align*} and equality is attained when $t_i=0$ for $1<i<n$.
\end{prop}

\begin{proof}
From \eqref{i}, we have 
\begin{align}\label{star}
& f(t_1,\dots, t_n) = t_1+2t_2+\dots + (n-1)t_{n-1} + n(1-t_1-\dots-t_{n-1})=\\
&\nonumber n-(n-1)t_1-(n-2)t_2-\dots-t_{n-1}
\end{align}

From \eqref{i} and \eqref{ii} we have
\[nt_1+(n-2)t_2+(n-3)t_3+\dots +t_{n-1}\ge n-1\]
which implies
\begin{align}
&t_1\ge \frac{n-1-(n-2)t_2-(n-3)t_3-\dots t_{n-1}}{n} \label{t1}
\end{align}

Applying \eqref{t1} to \eqref{star}, we get

\begin{align*}
& f(t_1,\dots, t_n) \le n - \frac{n-1}{n}[n-1-(n-2)t_2-(n-3)t_3-\dots-t_{n-1}]\\
& - (n-2)t_2-\dots - t_{n-1}=\\
&n-\frac{(n-1)^2}{n}+\Big{(}\frac{(n-1)(n-2)}{n}-(n-2)\Big{)}t_2+\\
&\Big{(}\frac{(n-1)(n-3)}{n}-(n-3)\Big{)}t_3+\dots +\Big{(}\frac{n-1}{n}-1\Big{)}t_{n-1}=\\
&\frac{2n-1}{n}t_1-\frac{n-2}{n}t_2-\frac{n-3}{n}t_3-\dots - \frac{1}{n}t_{n-1}\le \frac{2n-1}{n}t_1
\end{align*}
with equality when $t_2=\dots=t_{n-1}=0$.
\end{proof}

\section{A New Bound}

Our argument, in the spirit of Bartsalkin and German \cite{BG}, relies on labeling the vertices of a minimum dominating set, $D$, of $G\cart H$ with labels that contain integers from $\{1,\dots, \gamma(G)\}$. Labels may be sets of integers of size between one and $\pi(G)$. We show that every set of labels containing a fixed integer is at least of size $\gamma(H)$. We then control the overcount of vertices by applying the method of Contractor and Krop \cite{CK}. This is done by first applying a series of labeling refinements of the vertices of $D$. Labels may contain one or more integers and in each successive labeling, we attempt to reduce the number of labels with more than one integer while at the same time maintaining the property that vertices with labels that contain a fixed integer, when projected onto $H$, form a dominating set of $H$. To show a lower bound on $D$ we argue that every $G$-fiber with a certain number of vertices labeled by more than one integer must contain a fixed number of vertices labeled by one integer.

\begin{thm}\label{general}
For any graphs $G$ and $H$, 
\begin{align}
\gamma(G\square H)\geq \frac{\pi(G)}{2\pi(G)-1}\gamma(G)\gamma(H).\label{main}
\end{align}
\end{thm}

\begin{proof}
For any graphs $G$ and $H$, let $\Gamma=\{v_1,\dots, v_k\}$ be a minimum dominating set of $G$ and $D$ be a minimum dominating set of $G \cart H$.

Our proof is composed of increasingly refining labelings of the vertices of $D$, which successively reduce the number of labels with more than one element so that for any fixed label $i$, if we project vertices that contain $i$ in their label onto $H$, we produce a dominating set of $H$. We then show a bound on the label overcount to produce the desired inequality.

We now address the following property used in some portions of our labeling scheme: for any maximal $S \subseteq [k]$, if $v^h\in D\cap P^h_{S}$, then $v^h$ may only be labeled by a subset of $S$. For example, if $v^h\in P^h_{i,j}$, then $v^h$ may be labeled by $i,j,$ or $\{i,j\}$. We call labelings that follow this property \emph{faithful}. Note that not all of our labelings are faithful.

For any $h\in V(H)$, suppose the fiber $G^h$ contains $\ell_h(=\ell)$ vertically undominated cells $\big\{Q_{i_1}^h,\dots, Q_{i_{\ell}}^h\big\}$ for $0\leq \ell \leq k$. We set $I^h=\{i_1,\dots, i_{\ell}\}$ and $A^h=[k]-I^h$.
\medskip

We apply the procedure \emph{Labeling 1} to the vertices of $D$.\\

\underline{\bf Labeling $1$:}
\medskip

\begin{enumerate}[label=(\roman*)]
\item For $v^h\in D^h \cap Q_i^h$ for any $h\in V(H)$ and $1\leq i \leq k$, we label $v^h$ by $\{j: v^h \in N[v^h_j], j\in [k] \}$. Note that if $v^h\in P^h_i$, then the label assigned to $v^h$ is $\{i\}$. 
\item If $v^h\in D^h$ is a shared neighbor of some subset of $\{v^h_i:i\in I^h\}$, then it is a member of $P^h_S$ for some maximal $S\subseteq  I^h$, and we label $v^h$ by $S$. Note that in this designation, $v^h$ may also be a shared neighbor of some subset $Y$ of $A^h$.
\item If $v^h\in D^h$ is a member of $P^h_S$ where $S\subseteq A^h$, and $v^h$ is neither in $Q^h_i$ for $1\leq i \leq k$, nor a shared neighbor of any subset of $\{v^h_i:i\in I^h\}$, then do not label $v^h$ and call it \emph{available}. 
\end{enumerate}
This completes Labeling $1$.
\medskip

We relabel the vertices of $D$, doing so in $D^h$ for fixed $h\in V(H)$, stepwise, until we exhaust every $h\in V(H)$. This procedure, which we call \emph{Labeling 2}, is described next.\\

\underline{\bf Labeling $2$:}
\medskip

For every $h\in V(H)$, we list the labels of vertices of $D^h$, and write them in row $h$ in arbitrary order. This produces a two-dimensional array of $|V(H)|$ rows of labels, some of which may be empty. For an arbitrary $h\in V(H)$, we perform two alterations to the labels in row $h$ which we call the \emph{internal} and \emph{external alterations}. The purpose of the procedure is to reduce the number of vertices of $D$ with multiple integers in their labels. In each of these procedures, we note whether we follow \emph{the dominion rule}: if $v_i^h\in D^h$ with label $S$, then any alteration of $S$ must retain the label $i$. 

\medskip

Choose any pair of labels $S$ and $T$ in row $h$ and perform the \emph{internal alteration},
\begin{enumerate}
\item If $|S\cap T|>1$, then remove one common element from $S$ and another from $T$, arbitrarily, subject to the dominion rule. Repeat this step.
\item If $|S|=1$, $|T|>1$, and $S\cap T \neq \emptyset$, then remove the elements of $S$ from $T$, unless such an alteration defies the dominion rule. In that case, remove the element from $S$ and replace it with another element of $T$, which is removed from $T$, subject to the dominion rule. Note that if this case occurs, the labeling is not faithful.
\item If $|S|=|T|=1$, then make no changes to $S$ or $T$. 
\item If $|S|>1,|T|>1$, $|S\cap T|=1$, then remove the common element from one of $S$ or $T$, subject to the dominion rule. 
\end{enumerate}

We repeat this internal alteration for every row $h\in V(H)$ until every pair of labels in a row is mutually disjoint.
\medskip

We perform the \emph{external alteration} to the array obtained from the internal alteration. Choose any $h\in V(H)$ and suppose $N_H(h)=\{h_1,\dots, h_n\}$. For every label $S$ in row $h$, we consider labels $T$ of row $h_i$ for $i=1,\dots, n$, and repeat the relabeling from the internal alteration as follows:
\begin{enumerate}
\item Set $i=1$.
\item For every label $S$ in row $h$ and $T$ in row $h_i$, if $|S\cap T|>1$, then remove one common element from $S$ and another from $T$, arbitrarily, subject to the dominion rule. Repeat this step.
\item If $|S|=1$, $|T|>1$, and $S\cap T \neq \emptyset$, then remove the elements of $S$ from $T$, unless such an alteration defies the dominion rule. In that case, remove the element from $S$ and replace it with another element of $T$, which is removed from $T$, subject to the dominion rule. Note that if this case occurs, the labeling is not faithful.
\item If $|S|>1$, $|T|=1$, and $S\cap T \neq \emptyset$, then remove the label of $T$ from $S$, subject to the dominion rule.
\item If $|S|=|T|=1$, then make no changes to $S$ or $T$.
\item If $|S|>1,|T|>1$, and $|S\cap T|=1$, then remove the common element from one of $S$ or $T$ arbitrarily, subject to the dominion rule. 
\item Let $i=i+1$ and repeat this relabeling until $i=n+1$.
\end{enumerate}

After all alterations are performed for every row $h\in V(H)$, we confer the labels in the rows to the corresponding vertices of $D$. This completes Labeling $2$.

\medskip

%
%

\medskip

For $h\in V(H)$, let $R^h$ be the vertices of $D^h$ which have labels with more than one element. That is, vertices $v^h$ labeled by some $X$ so that $|X|>1$. Say $R^h=\{x_1^h,\dots, x_r^h\}$ where $|R^h|=r^h$. Let $m_i^h$ be the cardinality of the label on $x^h_i\in R^h$. For each vertex in $R^h$, we place each element from the label on that vertex in the set $J^h$. For example, if $R^h$ contains vertices with labels $\{i_1,i_2\}$ and $\{i_3,i_4\}$, then $J^h=\{i_1,i_2,i_3,i_4\}$. We call this process the \emph{label union}. 

\medskip

\begin{claim}\label{noshared}
For any $h\in V(H)$, $\C_{J^h}$ is dominated by vertices of $D^h$.
\end{claim}

\begin{proof}
Suppose $v^h\in \C_{J^h}$ is dominated by $u^{h'}\in D^{h'}$ for some $h'\in N_H(h)$. Note that since $v^h$ is in $\C_{J^h}$, the subset of vertices $v_1^h, \dots, v_k^h$, which are adjacent to $v^h$ may be only those which have an index in $J^h$. Furthermore, $u^{h'}$ is in the same chamber in $G^{h'}$ as $v^h$ is in $G^h$ since $v^h$ and $u^{h'}$ have the same $G$-coordinate. More precisely, projecting $v^h$ onto $G$ or $u^{h'}$ onto $G$, produces the same vertex $v\in V(G)$ such that $v\in \C_J$. Let $A$ be the label on $u^{h'}$. Since Labeling $2$ has been performed, for any vertex $w^h\in R^h$ with labels $T$, $A\cap T=\emptyset$ unless $|A|=|T|=1$. If $|T|=1$, then $w^h\notin R^h$, contradicting our assumption. Otherwise, elements of $A$ are not in $J^h$ leading to $v^h\notin \C_{J^h}$ which is a contradiction.
\end{proof}

%
%

\medskip

We now take advantage of all ``superfluous'' dominating vertices, which belong to the same cells and share the same singleton label. Let $S$ be a subset of vertices in $D^h$ all of which have singleton label $i$ for some $i\in[k]$. We call all but one of the elements of $S$ \emph{free}, choosing free vertices arbitrarily. We remind the reader that \emph{available} vertices were defined in part (iii) of Labeling $1$. We call the following procedure \emph{Labeling 3}.\\

\underline{\bf Labeling $3$:}
\medskip

For any set $F$ of free or available vertices in a given $G$-fiber, we relabel the vertices of $F$, whenever possible, by distinct labels of $J^h$, leaving at least one label in $J^h$, and remove those labels that had been used in the relabeling from the vertices of $R^h$. Call the sets thus produced $\hat{R^h}$ and $\hat{J^h}$, respectively. We then update the sets $R^h$ and $J^h$ so that $R^h=\hat{R^h}$ and $J^h=\hat{J^h}$. Vertices that were free or available and were relabeled this way, are no longer considered free or available. We iterate this procedure over all $h\in H$. This completes Labeling $3$.

\medskip

For convenience of notation, for any $v\in D$, let $\phi(v)$ be the label on $v$ after Labeling $3$. Furthermore, for any set of vertices $S\in D$, we write $\phi(S)$ to mean the set of labels on the vertices of $S$.

Set $B^h=D^h-(D^h\cap \C_{J^h})$. By Claim \ref{noshared} we can let $E_{J^h}$ be a minimum subset of vertices of $B^h$ so that $(D\cap \C_{J^h})\cup E_{J^h}$ dominates $\C_{J^h}$. That is, $E_{J^h}$ is a set of vertices with neighbors in $\C_{J^h}$, which along with vertices in $D\cap \C_{J^h}$, dominate $\C_{J^h}$. Note that by this definition, all vertices in $E_{J^h}$ must contain a singleton label and cannot be free or available. 

Note that any vertex $x^h$ of $E_{J^h}$ may be a private or shared neighbor of $\Gamma^h$ but may not itself be a member of $\Gamma^h$ since it dominates vertices in $\C_{J^h}$.

Furthermore, note that at this point, only vertices of $G^h-\C_{J^h}$ may be vertically dominated.

\begin{claim}\label{big}
For every $h\in V(H)$, $|E_{J^h}|\geq \sum_{i=1}^{r^h}(m_i^h-1)$.
\end{claim}

\begin{proof}
Set $j=|E_{J^h}|$ and notice that $E_{J^h}\cup R^h$ dominates $\C_{J^h}$. If we let $I'=[k]-J^h$, then $R^h\cup E_{J^h}\cup  (\bigcup_{i\in I'}v_i^h)$ dominates $G^h$. 
This yields

\begin{align*}
&|R^h\cup E_{J^h}\cup \Big{(}\bigcup_{i\in I'}v_i^h\Big{)}|= |R^h| + |E_{J^h}| + |\Big{(}\bigcup_{i\in I'}v_i^h\Big{)}| - |E_{J^h}\cap \Big{(}\bigcup_{i\in I'}v_i^h\Big{)}| \ge k
\end{align*}

Which is equivalent to

\begin{align*}
&r^h+j+k-\sum_{i=1}^{r^h} m_i^h- |E_{J^h} \cap \Big{(}\bigcup_{i\in I'}v_i^h\Big{)}| \ge k
\end{align*}

and produces

\begin{align*}
&\sum_{i=1}^{r^h} m_i^h\le r+j-|E_{J^h} \cap \Big{(}\bigcup_{i\in I'}v_i^h\Big{)}| \le r^h+j
\end{align*}

which is the desired inequality.
\end{proof}

Note that as a consequence of Claim \ref{big}, for any vertex $x_i^h\in R^h$ we can associate $m_i-1$ vertices of $E_{J^h}$ with $x_i$ such that for any $j\neq i$, the vertices associated to $x_i^h$ and $x_j^h$ are distinct.

\medskip

\begin{claim}\label{projection}
For a fixed $i$, $1\leq i \leq k$, projecting all vertices with labels containing $i$ to $H$ produces a dominating set of $H$.
\end{claim}

\begin{proof}
Let $S$ be the set of vertices of $H$ which are the projections of vertices with labels containing $i$. If $S$ is not a dominating set of $H$, then there exists a vertex $h\in V(H)$ which is not dominated by $S$. This means that 

\begin{enumerate}
\item $v_i^h\notin D$ since by the dominion rule, $i$ would be a label of $v_i^h$. 
\item $Q_i^h$ is vertically undominated, since otherwise $h$ would be dominated by $S$.
\item For any set $A\subseteq [k]$ such that $i\in A$, $P_A^h\cap D=\emptyset$ and $P_A^{h'}\cap D=\emptyset$ for all $h'\in N(h)$, since otherwise, Labelings $1$, $2$, and $3$ would ensure that either the label of some vertex in $P_A^h\cap D$ contains $i$ or the label of some vertex in $P_A^{h'}\cap D$ contains $i$ for some $h'\in N_H(h)$. In either case, $h$ would be dominated by $S$ in $H$.
\end{enumerate}

However, facts $(1),\, (2),$ and $(3)$ imply that $v_i^h$ is not dominated by $D$ in $G\cart H$, which is a contradiction to the assumption that $D$ is a dominating set.

\end{proof}

Call the set of vertices of $D$ containing the label $i$, $D_i$. Summing over all $i$ we count vertices with labels of cardinality $m_i^h$, $m_i^h$ times, for every $h\in V(H)$.

For any $j\in [k]$, let $F_j$ be the set of those vertices of $D$ with labels of cardinality $j$. We see that

\begin{align}
\gamma(G)\gamma(H) \leq \sum_{i=1}^k|D_i|=\sum_{j=1}^k j\times |F_j| =\sum_{j=1}^{\pi(G)} j\times |F_j| \label{eq1}
\end{align}


Define $t_i=\frac{|F_i|}{|D|}$, for $1\leq i \leq \pi(G)$. Note that $\sum_{i=1}^{\pi(G)}t_i=1$, which satisfies \eqref{i}. Furthermore, by Claim \ref{big}, we have

\begin{align*}
|F_1|\ge \sum_{h\in V(H)}|E_{J^h}| \ge \sum_{h\in V(H)}\sum_{i=1}^{r_h}(m_i^h-1)=\sum_{j=2}^{\pi(G)}(j-1)|F_j|
\end{align*}
which means that \eqref{ii} is satisfied as well.

We apply Proposition \ref{ineq} to find that  $(\ref{eq1}) \leq \frac{2\pi(G)-1}{\pi(G)}|D|$, and thus,

\[|D|\geq \frac{\pi(G)}{2\pi(G)-1}\gamma(G)\gamma(H).\]

\end{proof}

\subsection{Some Consequences}

By definition, $\pi (G)\leq \gamma(G)$ and $\pi(G)\leq \Delta(G)$, which can be trivially applied to \eqref{main}. However, we can improve these bounds slightly.

\begin{thm}\label{gammabound}
For any graphs $G$ and $H$,
\begin{align*}
&\gamma(G\cart H)\geq \frac{\lceil \frac{\gamma (G)}{2}\rceil}{2\lceil \frac{\gamma (G)}{2}\rceil-1}\gamma(G)\gamma(H),\, \text{for odd $\gamma(G)$} \text{ and}\\
&\gamma(G\cart H)\geq \frac{\gamma (G)}{2\gamma (G)-2}\gamma(G)\gamma(H),\, \text{for even $\gamma(G)$.}
\end{align*}
\end{thm}

\begin{proof}
We show that after the labelings in the proof of Theorem \ref{general}, for all $h\in H$, if $v^h\in R^h$, then the label on $v^h$ has cardinality at most $\lceil \frac{\gamma (G)}{2}\rceil$. For a contradiction, suppose that for some $h\in V(H)$, $v^h\in R^h$ has a label of cardinality at least $\lceil \frac{\gamma(G)}{2} \rceil+1$. By Claim \ref{big}, for any $h\in V(H)$, $|E_{J^h}|\geq \sum_{i=1}^{r}(m_i-1)$.


Claim \ref{big} allowed us to say that for any vertex $x_i^h\in R^h$ we can associate $m_i-1$ vertices of $E_{J^h}$ with $x_i$ such that for any $j\neq i$, the vertices associated to $x_i^h$ and $x_j^h$ are distinct. To be explicit, for any $h\in V(H)$, let $f^h: R^h \rightarrow E_{J^h}$ be the injective and surjective relation (which is not necessarily a function) that follows this association.

For any $u\in D^h$, if $\phi(u)$ is a singleton, we say that the cell $Q^h_{\phi(u)}$ is \emph{linked} if $v^h_{\phi(u)}\notin D^h$, and there exists $h'\in N_H(h)$ and a vertex $w\in D^{h'}$ so that $\phi(w)=\phi(u)$. We say that a cell $Q^h_{\phi(u)}$ is \emph{flat} if $v^h_{\phi(u)}\notin D^h$ and it is not linked. We say that a flat or linked cell $Q^h_i$ \emph{is guarded by} a vertex $v^h$ if $v^h$ is labeled by $\{i\}$. Let $C_i^h$ be the union of $Q_i^h$ and the vertices that guard $Q_i^h$, and call them the \emph{guarded cells} of $G^h$. 

\medskip

\noindent\emph{Comment $1$: Any flat cell $Q^h_i$ must be dominated by $D^h$ and cannot contain vertices of $\Gamma^h\cap D^h$.}

\medskip

\noindent\emph{Comment $2$: To motivate the definition of guarded cells, note that after Labeling $3$, there may be vertices of $D^h$ which are not adjacent to $v_i^h$ but are labeled $i$. We would like to consider such vertices, which would be guards of $Q_i^h$, as part of that cell, and hence define the guarded cells.}

\medskip

For a fixed $h\in V(H)$, call the set of linked cells in $G^h$, $L^h$, and the set of flat cells $F^h$. For any cell $C^h_i$, we write $||C^h_i||$ to mean the number of vertices in $D^h$ with label $\{i\}$, that is $|\{v^h:\phi(v^h)=i\}|$, and call this amount the \emph{weight} of $C^h_i$. In general, we define 

\begin{align*}
||L^h||=\sum_{C^h_i \text{linked}}||C^h_i|| \qquad \text{ and } \qquad ||F^h||=\sum_{C^h_i \text{flat}}||C^h_i||.
\end{align*}

Note here that $||L^h||$ and $||F^h||$ need not be the same as $|L^h|$ and $|F^h|$, respectively, since the latter refers to the number of linked or flat cells, not their weights.

\begin{claim}\label{flat}
\[(E_{J^h}\cup R^h) \cup [\Gamma^h\backslash (\{v^h_j:j\in J^h\}\cup \{v_i^h:C_i^h \in F^h\})]\] is a dominating set of $G^h$.
\end{claim}

\begin{proof}
By Claim \ref{noshared}, $E_{J^h}\cup R^h$ dominates $\C_{J^h}$ and any flat cell is dominated in $G^h$ by definition. We must now consider other vertices in the chambers of multiple flat cells, the chambers of flat and linked cells, and other chambers we have not yet considered. 

For $v^h\in \C_S$ where $S$ is some set of indices of $F^h$, suppose $v^h$ is not dominated by $D^h$. Then there must exist some $h'\in N_H(h)$ and a vertex $v^{h'}\in D^{h'}$ so that $v^h$ and $v^{h'}$ have the same $G$-coordinate. Note that during the external alteration of Labeling $2$, $v^{h'}$ was relabeled and received a label with elements from the set of indices of $S$. Suppose this label is $S'$. Since Labeling $2$ has been performed, either $S\cap S'=\emptyset$ or $|S|=|S'|=1$. This means that we may only consider the latter case, which implies that for some $i\in S'$, $C^h_i$ is a linked cell which contradicts our assumption that $C^h_i$ is flat. Hence, $v^h$ must be dominated by $D^h$. 

Next, if $v^h\in \C_T$ where $T$ contains a mixture of indices from $J^h$ and $F^h$, and $v^h$ is not dominated by $D^h$, then again there must exist some $h'\in N_H(h)$ and a vertex $v^{h'}\in D^{h'}$ so that $v^h$ and $v^{h'}$ have the same $G$-coordinate. Again, during the external alteration of Labeling $2$, $v^{h'}$ was relabeled to a label with elements from the set of indices of $T$. That label could not have contained elements from $J^h$ since the labels in $J^h$ are the remaining labels on the vertices of $R^h$ that were not singletons after the external alteration, and if $v^{h'}$ held such a label, then those elements could not have remained in $J^h$. Also, such a label could not have been in $F^h$ as in the previous case. Hence, $v^h$ must be dominated by $D^h$. 

Finally, suppose $v^h$ is in $\C_U$ where $U$ is an index set with no indices from the set $J^h$ or of $F^h$. In other words, $v^h$ is not in the chambers considered above. If $v^h$ is in a cell, say $C^h_j$, and this cell is flat, then $v_j^h\notin D$, since flat cells do not contain vertices of $\Gamma\cap D$. Thus, $v^h$ is dominated by $\Gamma^h\backslash (\{v^h_j:j\in J^h\}\cup \{v_i^h:C_i^h\in F^h\})$. If $v^h$ is in cell $C^h_j$ and this cell is not flat, then $v^h_j$ is in $\Gamma^h\backslash (\{v^h_j:j\in J^h\}\cup \{v_i^h:C_i^h\in F^h\})$ and $v^h_j$ dominates $v^h$.
\end{proof}

Claim \ref{flat} immediately yields the inequality

\begin{align*}
\gamma(G)&\le \big{|}(E_{J^h}\cup R^h) \cup [\Gamma^h\backslash (\{v^h_j:j\in J^h\}\cup \{v_i^h:C_i^h \in F^h\})]\big{|}\\
&=|(E_{J^h}\cup R^h)| + |\Gamma^h| - |\{v^h_j:j\in J^h\}| - |F^h|\\
&\le ||L^h||+||F^h||+|R^h|+\gamma(G)-|J^h|-|F^h|\\
&=||L^h||+||F^h|| - |F^h| +\gamma(G) - \sum_{i=1}^{r^h}(m_i^h-1)
\end{align*}

which gives us the inequality

\begin{align}
||L^h||+||F^h||-|F^h|\ge  \sum_{i=1}^{r^h}(m_i^h-1) \label{flatdom}
\end{align}

\medskip

Note that for any $h\in V(H)$ and fiber $G^h$, if $|R^h|>0$, then no flat cell of $F^h$ can have more than one vertex of $D$, since one of these vertices would have been free and relabeled by Labeling $3$. This observation implies that for all such fibers, $||F^h||-|F^h|=0$ and applying this to \eqref{flatdom}, we have

\begin{align}
||L^h|| \ge \sum_{i=1}^{r^h}(m_i^h-1) \label{flatdomer}
\end{align}

Thus, we may reassign vertices of $L^h$ to $E^h$, with at least $\sum_{i=1}^{r^h}(m_i^h-1)$ assigned. We perform this reassignment for every $h\in V(H)$.

\medskip

We now apply a procedure in which we relabel all vertices in $R^h$ by singleton labels, and increase the number of labels on vertices in a subset of $f(R^h)$ by one. In the process, we maintain the property that projections of vertices of $D$ containing a fixed label onto $H$ produce a dominating set of $H$. Inequality \eqref{flatdom} ensures that there is enough weight in the flat cells of any fiber so that all but one of the dominating vertices may be relabeled by labels of $R^h$, which subsequently would be removed from $R^h$, until the remaining labels on $R^h$ are singletons. This is done after the dominating vertices of linked cells each receive one additional label from $R^h$, which is subsequently removed from $R^h$.

\medskip

We call this procedure \emph{Labeling 4}.

\medskip
\underline{\bf Labeling $4$:}
\medskip

For any $h\in V(H)$, every vertex $x_i^h \in R^h$ and associated vertices from linked cells in $f(x_i^h)\cap L^h \subseteq E_{J^h}$, add to the label of each vertex $v\in f(x_i^h)\cap L^h$ a distinct element of $\phi(x_i^h)$ and remove that label from $x_i^h$. Note that the relabeled vertices of $f(x_i^h)\cap L^h$ will each have a label with two elements.

This completes Labeling $4$.

\medskip

Notice that after Labeling $4$, the only vertices of $D$ that may not have a singleton label are the dominating vertices of linked cells. Furthermore, all labels that are not singletons contain exactly two elements. 

\medskip

For $i\in [k]$, let $L_i$ be the set of linked cells of $\bigcup_{h\in V(H)}\Big{(}\bigcup_{C^h\in L^h}C^h\Big{)}$, which  contain the label $i$ after Labeling $4$. Notice that the vertices of $D$ in these cells are the only vertices that do not have a singleton label. Let $W_i$ be the set of vertically dominated cells from $\bigcup_{h\in V(H)}Q^h_i$ which contain no vertices of $D$ and which, when projected onto $G$ contain $v_i$. Let $V_i$ be the set of cells from $\bigcup_{h\in V(H)}Q^h_i$ so that every cell contains an element from $\Gamma^h$ for some $h\in V(H)$.

\medskip


Furthermore, note that after Labeling $4$, for any $h\in V(H)$, $G^h$ cannot contain both a vertex of $R^h$ and a free vertex. This implies that for each $i\in I-J^h$, $D^h$ may contain at most one vertex with label $i$. However, note that $v^h$ is adjacent to $\lceil \frac{\gamma(G)}{2} \rceil+1$ vertices of $v_1^h,\dots, v_k^h$ which have different labels from the vertices of $E_{J^h}$, which are themselves labeled distincly by at least $\lceil \frac{\gamma(G)}{2} \rceil$ singleton labels. Summing these distinct labels produces $2\lceil \frac{\gamma(G)}{2} \rceil+1$ which is larger than $\gamma(G)$, a contradiction. Thus, no label of any vertex in $D$ can be of cardinality greater than $\lceil \frac{\gamma(G)}{2} \rceil$. Applying the counting argument starting at \eqref{eq1} with $\pi(G)$ substituted by $\lceil \frac{\gamma(G)}{2} \rceil$ produces the desired conclusion.

The proof of the second inequality is identical.
\end{proof}

We note that this bound shows that Vizing's conjecture holds for graphs $G$ with domination number $2$ and any $H$, which was proved in \cite{BG}. Using the methods we developed in Theorem \ref{general}, it is not difficult to extend this result to the case when $\gamma(G)=3$ (which was solved in \cite{B}).


\begin{cor}\label{3}
For any graphs $G$ and $H$ so that $\gamma(G)=3$, 
\[\gamma(G\cart H)\geq \gamma(G)\gamma(H).\]
\end{cor}

\begin{proof}
First, we consider the case when $\gamma(G)=3$.
We perform the labelings from the proof of Theorem \ref{general}. By the proof of Theorem \ref{gammabound}, we have $\pi(G)\le \lceil\frac{\gamma (G)}{2}\rceil$. Since $\gamma(G)=3$, we have that $\pi(G)\le 2$. If all vertices have a singleton label, then we can show the proposed inequality by projecting all vertices with a fixed label onto $H$ as in the proof of Theorem \ref{main}. Hence, we assume that for some vertex $v^{h_1}\in D$ for $h_1\in V(H)$, has a paired label, and without loss of generality, say that label is $\{1,2\}$.
This means that $v^{h_1}$ is adjacent to $v_1^{h_1}$ and $v_2^{h_1}$. By Claim \ref{big}, $E_{J^{h_1}}$ must contain at least one vertex of $D$, $u^{h_1}$, with a singleton label, say $3$. Note that there can be no other vertices with label $3$ in $G^{h_1}$ since we relabeled all free and available vertices in Labeling $3$, after which there can be no $G$-fiber containing a free vertex and a vertex of $D$ labeled by more than one integer. Explicitly, $u^{h_1}$ and $v^{h_1}$ are uniquely associated by the function $f^{h_1}$. 

Note that $Q_3^{h_1}$ is a linked cell guarded by $u^{h_1}$ and say that it is linked to $Q_3^{h_2}$ which is guarded by $w^{h_2}$ for some $h_2\in V(H)$. 

Note that if $w^{h_2}$ were associated (by the function $f^{h_2}$) to a vertex of $D\cap G^{h_2}$ with a paired label, then that label would have to contain the pair of integers $\{1,2\}$. However, if either $1$ or $2$ were in that label, then Labeling $2$ would have reduced that label either in that vertex or in $v^{h_1}$, which is a contradiction. Thus, $w^{h_2}$ is not associated with any vertex of $D\cap G^{h_2}$. 

Project all vertices of $D$ that contain $3$ in their labels onto $H$. We now consider neighbors of $h_1$ in $H$ that are only dominated by $h_1$. Say $h_3$ is such a vertex of $H$. In this case, $Q_3^{h_3}$ is vertically dominated only by $u^{h_1}$.

Since $u^{h_1}$ does not dominate $v_3^{h_3}$, there must exist some vertex $y^{h_3}\in D$ which horizontally dominates $v_3^{h_3}$. However, since $h_1$ and $h_3$ are adjacent in $H$, $v^{h_1}$ is labeled by $\{1,2\}$ after Labeling $3$, and $u^{h_1}$ is labeled by $3$, we conclude that there can be no such $y^{h_3}$. Thus, we may relabel $u^{h_1}$ by $2$ and $v^{h_1}$ by $1$, noting that after this relabeling, the projection onto $H$ of all vertices with a $3$ in their label still dominates $H$.

By this observation, we may relabel all vertices with paired labels to singletons, and the result follows.

%

\end{proof}

 \bibliographystyle{plain}
 
 \end{document}